\documentclass[11pt]{amsart}

\newcommand{\N}{\mathbb{N}}
\newcommand{\Z}{\mathbb{Z}}
\newcommand{\R}{\mathbb{R}}

\def\di{\displaystyle}
\usepackage{amssymb}
\usepackage{latexsym}
\usepackage{xypic}
\usepackage{eufrak}
\usepackage{euscript}
\newtheorem{thm}{Theorem}[section]
\newtheorem{defi}{Definition}[section]

\newtheorem{cor}{Corollary}[section]

\newtheorem{rema}{Remark}[section]
\def\cal{\mathcal}
\begin{document}
\baselineskip 6mm
\title{About analytic non integrability}
\author{Jacky Cresson}
\address{Universit\'e de Franche-Comt\'e, \'Equipe de Math\'ematiques de Besan\c{c}on, CNRS-UMR 6623, 16 route de Gray, 25030
Besan\c{c}on cedex, France.}
\email{cresson@math.univ-fcomte.fr}
\maketitle

\begin{abstract}
We discuss non existence of analytic first integrals for analytic diffeomorphisms possessing a
hyperbolic fixed point with a homoclinic connection.
\end{abstract}

\tableofcontents

\section{Introduction}

Since the proof by J\"urgen Moser \cite{mo} of non existence of analytic first integrals for diffeomorphisms of the plane possessing a
hyperbolic fixed point with a transverse homoclinic crossing, many efforts have been devoted to generalize his approach.
The basic idea behind the Moser analytic non integrability theorem is that there exists a very complicated set (a Cantor set)
on which the dynamics is conjugated to a Bernoulli shift and this must prevent the system to possess an analytic first
integral. This idea goes back to Henri Poincar\'e \cite{po} in his famous work on the three body problem\footnote{``Que
l'on cherche \`a se repr\'esenter la figure form\'ee par ces deux courbes et leurs intersections en nombre infini dont chacune
correspond \`a une solution doublement asymptotique, ces intersections forment une sorte de treillis, de tissu, de reseau
\`a mailles infiniment serr\'ees; chacune de ces courbes ne doit jamais se recouper elle-m\^eme, mais elle doit se replier
elle-m\^eme d'une mani\`ere tr\`es complexe pour venir recouper une infinit\'e de fois toutes les mailles du r\'eseau. On sera
frapp\'e de la complexit\'e de cette figure, que je ne cherche m\^eme pas \`a tracer. Rien n'est plus propre \`a nous donner
une id\'ee de la complexit\'e du probl\`eme des 3 corps et en g\'en\'eral de tous les probl\`emes de la dynamique o\`u il
n'y a pas d'int\'egrale uniforme."}. We refer to the work of Dovbysh \cite{dov} for a continuation
of Moser's original approach.\\

However, as proved in \cite{cr}, we can deduce analytic non integrability without
any reference to this invariant Cantor set, only regarding on the local dynamics on the stable and unstable manifold and
a trivial argument using the transversality of the intersection. This proves that obstruction to analytic integrability is
related to phenomena which need not to be complex from the dynamical point of view. In that respect, I think that it is
misleading to base a theorem of analytic non integrability on the existence of an invariant hyperbolic set on which the
dynamics is conjugated to the shift map. The correct obstruction is not yet founded.\\

The main advantage of the method \cite{cr} is that its
generalization to more general invariant sets (normally hyperbolic compact manifolds, partially hyperbolic tori) is
easy. We refer to (\cite{cr},Theorem 6.1,p.299) for a general theorem concerning normally hyperbolic compact manifolds with
a transversal homoclinic intersection using our method. \\

In the case of hyperbolic point, our result can be stated as follow:

\begin{thm}[\cite{cr},Theorem 2.2,p.291]
Let $f$ be an analytic diffeomorphism of $\R^n$ possessing a hyperbolic fixed point $p$, with stable and unstable
manifolds denoted by $W^- (p)$ and $W^+ (p)$ respectively. We denote by $\lambda^{\sigma} =(\lambda_i^{\sigma})_{i=1,\dots ,
n^{\sigma}}$ the eigenvalues of $Df(p)$ associated to $W^{\sigma} (p)$, with
$n^{\sigma} =\dim W^{\sigma} (p)$, $\sigma=\pm$, $n^- +n^+ =n$. We assume that:\\

i) $W^- (p)$ and $W^+ (p)$ intersect transversally at a homoclinic point $h$,

ii) For $\sigma =\pm$, the spectrum $\lambda^{\sigma}$ satisfy a multiplicative non resonance condition, i.e.
\begin{equation}
\mid (\lambda^{\sigma})^{\nu} \mid \not= 1,\ \mbox{\rm for all}\ \nu \in \Z^{n^{\sigma}}\setminus \{ 0\} .
\end{equation}

iii) The homoclinic point is admissible\footnote{We refer to \cite{cr} for a definition. By ii), we can linearize analytically
$f$ in a neighbourhood of $p$ on $W^{\sigma} (p)$, $\sigma$. We denote by $U^{\sigma}$ such a neighbourhood of linearization. We denote
by $h^{\sigma}$ the first image of $h$ by a given iterate of $f$ in $U^{\sigma}$. The admissibility condition is equivalent to
impose that all the coordinates of $h^{\sigma}$ restricted to $W^{\sigma} (p)$ in the linearizing coordinates system are
non zero. A problem is to find a suitable geometric interpretation of this condition.}.\\

Then $f$ does not admit a non trivial analytic first integral.
\end{thm}

The transversality assumption is fundamental in the proof of this theorem proposed in \cite{cr}. However, we already know
some particular examples where this assumption can be relaxed. We refer to Cushman \cite{cu} for symplectic diffeomorphisms
of the plane.\\

The condition on the spectrum of $f$ can be relaxed. Indeed, in some resonant cases, we have generalized this result (see
\cite{cdr}). However, we do not know if we can cancel this condition completely. This is the case for $n=2$.\\

In this paper, we develop a new method to prove non existence of analytic first integrals, inspired by a previous work of
Cushman \cite{cu}, which allows us to relax the transversality and the resonance assumptions. We also discuss some classical
ideas about analytic non integrability and the positivity of topological entropy.

\section{Topological crossing and topological entropy}

In this section, we discuss several extension of our theorem which relax the transversality assumption, using a recent
result of Rayskin \cite{ray}. We also discuss an alternative formulation using topological entropy and discuss its
significance to classify analytic integrable and non integrable dynamical systems.

\subsection{Topological crossing and Rayskin's theorem}

We begin with some preliminary definitions and notations.

\subsubsection{Graph portion}

Let $f$ be a diffeomorphism of $\R^n$ with a hyperbolic fixed at the origin and $U\subset \R^n$ be some small neighbourhood
of $0$. Denote by $W^-$ (resp. $W^+$) the associated stable (resp. unstable) manifold, and by $n_-$ (resp. $n_+$) its
dimension ($n_- +n_+ =n$). Let $h$ be a homoclinic point of $W^-$ and $W^+$. Denote by $\cal V$ a small $n_+$-neighbourhood
in $W^+$ around the origin. Define a local coordinate system $E_+$ at $0$ which spans $\cal V$. Similarly, define a local
coordinate system $E_-$ which spans $W^-$ near $0$. Let $E=E_- +E_+ \subset U$.

\begin{defi}
A $n_+$-neighbourhood $\Lambda \subset W^+$ is a graph portion in $U$, associated with the homoclinic point $h$, if for some
$n>0$, $f^n (h)\in \Lambda$, $\Lambda \subset (U\cap W^+ )$, and $\Lambda$ is a graph in $E$-coordinates of some $C^1$-function
defined on $\cal V$.
\end{defi}

\subsubsection{Newhouse type}

Let $f$ be a diffeomorphism on $\R^n$ which has $1$-dimensional stable and $(n-1)$-dimensional unstable manifolds, and a
homoclinic point $h$. Since the stable manifold is $1$-dimensional, we can define one- and two-sided homoclinic
intersections.

\begin{defi}
We say that $f$ has a homoclinic tangency of Newhouse type, if $f$ possesses only one-sided homoclinic tangencies.
\end{defi}

\subsubsection{Main result}

The main result of Rayskin is the following (\cite{ray},Theorem 3.2,p.475):

\begin{thm}[Rayskin,2005]
Let $M$ be a smooth manifold, and let $f:M\rightarrow M$ be a diffeomorphism with a hyperbolic fixed point $p\in M$
and with stable and unstable invariant manifolds $W^-$ and $W^+$ at $p$. Assume:

(i) $f$ is locally $C^1$ linearizable in a neighbourhood of $p$;

(ii) $W^-$ and $W^+$ have a complementary dimension $n-1$ and $1$;

(iii) $W^-$ and $W^+$ have a finite tangential contact at an isolated homoclinic point $h\in M$ ($h\not= p$),
which is not of Newhouse type;

(iv) $W^+$ contains a graph portion $\Lambda$ associated with $h$ in a small neighbourhood, contained in the
neighbourhood of linearization.

Then the invariant manifolds have a point of transverse intersection arbitrary close to the point $h$.
\end{thm}

Assumption ii) is a priori not necessary as well as the condition to have only a finite order of contact. We refer to
\cite{cfp} for more details.

\subsection{Analytic first integrals and topological crossing}

In this first section, we extend our previous result on non existence of analytic first integrals for diffeomorphims
possessing a hyperbolic fixed point using the result of V. Rayskin \cite{ray}. The main point is that we are able to
replace the transversality assumption by a weaker one, precisely the existence of a topological crossing.

\begin{thm}
\label{maintopo}
Let $f$ be an analytic diffeomorphism of $\R^n$ such that $p$ is a hyperbolic fixed point for $f$, with stable and
unstable manifolds denoted by $W^- (p)$ and $W^+ (p)$ respectively which have complementary dimensions $n-1$ and $1$.
We denote by $\lambda^- =(\lambda_i^-)_{i=1,\dots ,n^-}$
and $\lambda^+ =(\lambda_i^+ )_{i=1,\dots ,n^+}$ the eigenvalues of $Df(p)$ associated to $W^- (p)$ and $W^+ (p)$, with
$n^{\sigma} =\dim W^{\sigma} (p)$, $\sigma=\pm$. We assume that:\\

i) $f$ is locally $C^1$ linearizable in a neighbourhood $U$ of $p$,

ii) $W^- (p)$ and $W^+ (p)$ have a finite tangential contact at an isolated homoclinic point $h$, $h\not= p$, which is not of
Newhouse type,

iii) $W^+ (p)$ contains a graph portion $\Lambda$ associated with $h$ in a small neighbourhood $V\subset U$,

iv) We denote by $h^{\sigma}$ the first intersection of $W^{\sigma} (p)$ with $W^{-\sigma} (p)\cap U$. We assume that
$h^{\sigma}$ is admissible, i.e. that the coordinates of $h^{\sigma}$ on $W^{\sigma} (p)$ via the linearizing map obtained
under assumption i) are all non zero;

v) the eigenvalues $\lambda^{\sigma}$, $\sigma=\pm$, satisfy a multiplicative non resonance condition, i.e.
\begin{equation}
\mid (\lambda^{\sigma} )^{\nu } \mid \not= 1 ,\ \mbox{\rm for all}\ \nu \in \Z^{n^{\sigma}} \setminus \{ 0\} ,
\end{equation}
where $\nu =(\nu_1 ,\dots ,\nu_{n^{\sigma}} )$, $(\lambda^{\sigma} )^{\nu} =(\lambda_1^{\sigma} )^{\nu_1} \dots
(\lambda_{n^{\sigma}}^{\sigma} )^{\nu_{n^{\sigma}}} $.\\

Then, the discrete dynamical system defined by $f$ does not possess a non trivial analytic first integral.
\end{thm}

\begin{proof}
The proof is elementary and reduces to prove that the assumptions of (\cite{cr},Theorem 2.2 p.291) are satisfied.\\

Assumptions i), ii), iii) imply that Rayskin's result (\cite{ray},theorem 3.2,p.475) can be applied. As a consequence,
there exists a transverse intersection of the stable and unstable manifold arbitrary close to the point $h$. We
denote by $h'$ this point.

As $h$ is admissible and the admissibility condition is stable under $C^0$ small perturbations, we conclude that the point
$h'$ can be chosen such that $h'$ is also admissible.

This fact and assumptions iv) and v) imply that the assumption of (\cite{cr},Theorem 2.2) are satisfied. As a consequence,
$f$ does not admit non trivial analytic first integrals.
\end{proof}

\subsection{Analytic integrals and topological entropy}

In the two dimensional case, our assumptions can be simplified. We have:

\begin{thm}
\label{main3}
Let $f$ be an analytic diffeomorphism of $\R^2$ such that $p$ is a hyperbolic fixed point for $f$, with one dimensional
stable and unstable manifolds denoted by $W^- (p)$ and $W^+ (p)$ respectively. We assume that:\\

i) $W^- (p)$ and $W^+ (p)$ intersect at an isolated homoclinic point $B$, $B\not= p$,
which is not of Newhouse type.\\

Then, the discrete dynamical system defined by $f$ does not possess a non trivial analytic first integral.
\end{thm}

We note that the fact to restrict our attention to analytic diffeomorphisms implies that if the stable and unstable manifolds
do not coincide and intersect, this intersection is necessarily of finite order. Moreover, by a theorem of Hartman \cite{ha},
a diffeomorphism of $\R^2$ can be $C^1$ linearized near a hyperbolic fixed point.\\

Remark that this theorem can of course be deduced from Moser's theorem (\cite{mo}) using the fact that i) implies the existence of a
transverse homoclinic point in every neighbourhood of the degenerate homoclinic point as showed by Charles Conley \cite{co}
(see \cite{cr2} for a proof).\\

A convenient formulation of the result can be obtained using the notion of {\it topological entropy} \cite{young}. The diffeomorphisms
considered in theorem \ref{main3} have positive topological entropy (see \cite{bw},\cite{hw}). In fact, we have in the two
dimensional case a general result linking positivity of the topological entropy and analytic non integrability.

\begin{thm}
Let $f$ be an analytic diffeomorphism of $\R^2$ with a strictly positive topological entropy, then $f$ does not possess a non
trivial analytic first integrals.
\end{thm}

\begin{proof}
The proof uses Katok's theorem (\cite{k1},\cite{k2}) which ensures the existence of a horseshoes for $f$ which carry
most of the topological entropy of $f$. By Moser's theorem \cite{mo}, the existence of a horseshoes implies analytic non
integrability.
\end{proof}

Many people tries to relate the positivity of the topological entropy with non integrability. However, as proved by
Bolsinov and Taimanov \cite{bt} the positivity of topological entropy as a criterion for non integrability
or chaos is not correct.

\section{Non existence of analytic first integrals: the linear case}

In this section, we relax or assumptions on the spectrum of the diffeomorphisms that we consider. The idea is to use not only
the dynamics on the stable and unstable manifold, but on an open neighbourhood of the fixed point. Doing this, we are
able to consider problems which are impossible via our previous approach.

\begin{thm}
\label{main}
Let $f$ be an analytic diffeomorphism of $\R^n$, which possesses a hyperbolic point $p$, with a stable manifold denoted
by $W^- (p)$ and an unstable manifold $W^+ (p)$ of dimension $n_-$ and $n_+$ such that $n_- +n_+ =n$. We denote by
$\lambda^-_i$, $i=1,\dots ,n_-$ and $\lambda^+_i$, $i=1,\dots ,n_+$ the eigenvalues of $Df(p)$ associated to $W^- (p)$ and
$W^+ (p)$. We assume that:\\

i) $W^+ (p)$ and $W^- (p)$ do not coincide and intersect along (at least) one homoclinic orbit $\gamma$,

ii) There exists an analytic coordinates system $(x,y)\in \R^{n_-} \times \R^{n_+}$ which linearizes $f$ in a
neighbourhood $U$ of $p$.

iii) We denote by $h^{\sigma}$ the first intersection of $\gamma$ with $W^{\sigma}\cap U$. We have $h^- =(x_- ,0)$ and
$h^+ =(0,y_+ )$. The point $h^-$ (resp. $h^+$) is admissible if all the coordinates of $x_-$ (resp. $y_+$) are non zero.

We assume that the spectrum $(\lambda^{\sigma})$ satisfies a multiplicative non resonance condition and $h_i^{\sigma}$ is
admissible, for $\sigma =+$ or $-$.

iv) $W^+ (p)$ contains a graph portion $\Lambda$ associated with $h$ in a small neighbourhood $V\subset U$.\\

Then, the dynamical system defined by $f$ does not possess a non trivial analytic first integral.
\end{thm}

The main differences with theorem \ref{maintopo} are the following:\\

- First, we do not assume that the multiplicative non resonance condition applies on the stable and unstable manifold but only
on one of it.\\

- Second, we do not assume a topological crossing, which is a very particular case of intersection between the stable and
unstable manifold, nor that one of the stable or unstable manifold is a curve. We can for example assume that $f$
has a homoclinic tangency of Newhouse type (see \cite{ray},definition 3.1). In that case, Rayskin's result on the
existence of a transversal intersection arbitrary close to the degenerate intersection does not apply and we can not use
the method developped in \cite{cr}.\\

The proof of theorem \ref{main} is inspired by the Cushman proof \cite{cu} of non existence
of analytic first integral for symplectic diffeomorphism of the plane possessing a finite contact homoclinic point and is
given in the next section.\\

We have the following comment on assumption ii) and iv):\\

- Assumption ii) is of technical nature. It simplifies the proof but is not necessary. We explain more precisely the
problem that we must solve in the next section.\\

- Assumption iv) is not a priori necessary. It is related to a construction which is of topological nature (see $\S$.
\ref{preliminary}) in the proof of theorem \ref{main}.\\

An interesting corollary of theorem \ref{main} concerns the case where either the stable or unstable manifold is one
dimensional.

\begin{cor}
Let $f$ be an analytic diffeomorphism of $\R^n$, which possesses a hyperbolic point $p$, with a stable manifold denoted
by $W^- (p)$ and an unstable manifold $W^+ (p)$ of dimension $n-1$ and $1$. We denote by
$\lambda^-_i$, $i=1,\dots ,n-1$ and $\lambda^+$, the eigenvalues of $Df(p)$ associated to $W^- (p)$ and
$W^+ (p)$. We assume that:\\

i) $W^- (p)$ and $W^+ (p)$ do not coincide and intersect along (at least) one homoclinic orbit $\gamma$,

ii) There exists an analytic coordinates system $(x,y)\in \R^{n -1} \times \R$ which linearizes $f$ in a
neighbourhood $U$ of $p$.

iii) $W^+ (p)$ contains a graph portion $\Lambda$ associated with $h$ in a small neighbourhood $V\subset U$.\\

Then, the dynamical system defined by $f$ does not possess a non trivial analytic first integral.
\end{cor}

Indeed, in that case, the spectrum on the one dimensional stable or unstable manifold satisfies the multiplicative non
resonance condition the homoclinic point is necessarily admissible.

\section{Proof of theorem \ref{main}}

We write the proof assuming that the spectrum on the stable manifold satisfies the multiplicative non resonance condition. The other
case is easily deduced using $f^{-1}$ instead of $f$.

\subsection{Preliminaries}
\label{preliminary}

In $U$, $f$ is analytically conjugated to
\begin{equation}
f(x,y) =(\lambda^- x ,\lambda^+ y ) , \ \ (x,y)\in \R^{n^-} \times \R^{n^+} .
\end{equation}
Let $C_+=(0,Y_+ )$, $Y_+ \not = 0$. We denote by $P_{C_+}$ the affine space defined by
\begin{equation}
P_{C_+} =C_+ +E^- ,
\end{equation}
where $E^-$ is the vectorial space $\R^{n^-} \times \{ 0\}$.\\

We denote by $h$ the first intersection of $W^+ (0)$ and $W^- (0)$ in $U$. We have $h=(x_- ,0)$ with $x_- \not= 0$. We can
always find a sequence of points $M_i =(x_i ,y_i)$ in $W^+ (0)$ such that
\begin{equation}
f^i (M_i ) \in P_{C_+} ,
\end{equation}
for $i$ greater than a given constant $N$ depending on $C_+$.\\

Indeed, as $W^+ (0)$ contains a graph portion $\Lambda$, we have up to choose $U$ sufficiently small a parametrization of
$\Lambda$ which is given by
\begin{equation}
y\in U_+ \subset \R^{n^+},\ \ y\longmapsto (\Lambda (y) ,y) ,
\end{equation}
where $U^+$ is an open neighbourhood of $0$. By definition, we have
\begin{equation}
\Lambda (0)=x_- ,
\end{equation}
but we have no assumptions on the derivative of $\Lambda$ at $0$. As a consequence, we have no assumption on the nature of
the intersection (transversal or not).\\

Let $(x_i ,y_i )$ such that $f^i (x_i ,y_i )\in P_{C_+}$. We then have
\begin{equation}
y_i =\lambda_+^{-i} Y_+\ \ \mbox{\rm and}\ \ x_i =\Lambda (\lambda_+^{-i} Y_+ ) .
\end{equation}
We have of course $x_i \rightarrow x_-$ when $i$ goes to infinity. We then have
\begin{equation}
X_i =\lambda_-^i \Lambda (\lambda_+^{-i} Y_+ ) .
\end{equation}
We have $X_i \rightarrow 0$ when $i$ goes to infinity.

\begin{rema}
The previous construction of a set of point on $P_{C_+}$ can certainly be done without any assumption about the existence of
a graph portion in $W^+ (0)$. This is essentially a topological phenomena which can be deduced from the Hartman-Grobman
theorem for hyperbolic points.
\end{rema}

We denote $N_i =f^i (M_i ) =(X_i ,Y_+)$.

\subsection{First step}

Let $I(x,y)$ be an analytic first integral for $f$. For $U$ sufficiently small, we have
\begin{equation}
I(x,y)=\di\sum_{(\nu,w )\in \N^{n^-} \times \N^{n^+} } a_{\nu ,w} x^{\nu} y^w ,
\end{equation}
with $a_{\nu ,w} \in \R$ for all $(\nu ,w)\in \N^{n^-} \times \N^{n^+}$.\\

As $N_i \in W^+ (0)$ for all $i$, we have $I(N_i )=I(h)=e$, where $e \in \R$ is a constant. We then have
\begin{equation}
\label{first}
I(N_i )= \di\sum_{(\nu,w )\in \N^{n^-} \times \N^{n^+} } a_{\nu ,w} X_i^{\nu} Y_+^w =e ,
\end{equation}
for all $i$ sufficiently great.\\

We have $X_i =\lambda^i x_i$ and $x_i \rightarrow x_-$ when $i$ goes to infinity. We assume that $\lambda$ satisfies the
multiplicative non resonance condition. Then, the quantities $\lambda^{\nu}$ can be totally ordered:
\begin{equation}
\lambda_-^{\nu_0 =0}=1 >\lambda_-^{\nu_1} > \lambda_-^{\nu_2} >\dots >\lambda_-^{\nu_j} >\dots .
\end{equation}
Equation (\ref{first}) can then be written as follow
\begin{equation}
\di\sum_{j\in \N} \lambda_-^{i\nu_j}  \left [
\di\sum_{w\in \N^{n_+}} a_{\nu_j ,w} x_i^{\nu_j} Y_+^w \right ] =e ,
\end{equation}
for all $i$ sufficiently great. As a consequence, we obtain by iteration the following set of equations
\begin{equation}
\label {second}
\left .
\begin{array}{l}
\di\sum_{w\in \N^{n_+}} a_{0,w} Y_+^w =e ,\\
\di\sum_{w\in \N^{n_+}} a_{\nu_j ,w} x_-^{\nu_j} Y_+^w =0 .
\end{array}
\right .
\end{equation}

\subsection{Second step}

Equation (\ref{second}) are valid for all $(0,Y_+ )\in W^+ (0)$. As a consequence, we can choose a sequence of
$Y_+$ of the following form
\begin{equation}
Y_+ (n)=\delta^n C ,
\end{equation}
where $\delta \in \R^{n_+}$, $0<\mu_i <1$ for $i=1,\dots ,n_+$, and $C$ is a constant vector such that $(0,C)\in W^+ (0)$,
$C_i \not=0$ for $i=1,\dots ,n_+$, and $\delta$ satisfies a multiplicative non resonance condition, i.e.
\begin{equation}
\delta^{\nu} \not=1 ,\ \ \forall \ \nu\in \Z^{n_+} \setminus \{ 0\} .
\end{equation}
For all $Y_+ (n)$, we then deduce from equation (\ref{second}):
\begin{equation}
\label{third}
\left .
\begin{array}{l}
\di\sum_{w\in \N^{n_+}} a_{0,w} \delta^{nw} C^w =e ,\\
\di\sum_{w\in \N^{n_+}} a_{\nu_j ,w} x_-^{\nu_j} \delta^{nw} C^w =0 .
\end{array}
\right .
\end{equation}
By the multiplicative non resonance assumption on $\delta$, we have a total order for $\delta^w$, $w\in \N^{n_+}$, i.e.
\begin{equation}
\delta^{w_0 =0} =1 >\delta^{w_1} >\dots >\delta^{w_j} >\dots .
\end{equation}
We then obtain
\begin{equation}
\label{fouth}
\left .
\begin{array}{l}
\di\sum_{k\in \N} a_{0,w_k } \delta^{nw_k } C^{w_k} =e ,\\
\di\sum_{k\in \N} a_{\nu_j ,w_k} x_-^{\nu_j} \delta^{nw_k} C^{w_k} =0 ,
\end{array}
\right .
\end{equation}
for all $j$ and $k$.
We deduce easily that
\begin{equation}
a_{0,0} =e,\ \ a_{\nu_j ,w_k } =0,\ \forall (j,k)\in \N \times \N \setminus (0,0) .
\end{equation}
As a consequence, the first integral is trivial, i.e. constant on $U$.

\section{Non existence of analytic first integrals: a general case}

The most stringent assumption of theorem \ref{main} concerns the analytic linearization in a neighbourhood of the hyperbolic
fixed point. However, as we will see, this condition
can be cancelled. This is already the case in the Cushman theorem \cite{cu}.

\subsection{Set-up of the problem}

Assume that $f$ has the following general form
\begin{equation}
f(x,y)=(\lambda_- x ,\lambda_+ y) +(r_- (x,y) ,r_+ (x,y)) .
\end{equation}
We denote by
\begin{equation}
f^n (x,y)=(\lambda_-^n x,\lambda_+^n y) +(r_{-,n}(x,y) ,r_{+,n} (x,y)) .
\end{equation}
As $(X_i ,Y_+) =f^i (x_i ,y_i )$, equation (\ref{first}) is then equivalent to
\begin{equation}
\label{exten}
I(N_i )= \di\sum_{(\nu,w )\in \N^{n^-} \times \N^{n^+} } a_{\nu ,w} (\lambda_-^i x_i +r_{-,i} (x_i ,y_i ) )^{\nu} Y_+^w =e ,
\end{equation}
for all $i$ sufficiently great.\\

If we assume that there exists a limit to
\begin{equation}
\label{condi}
\lambda_-^{-i} r_{-,i} (x_i ,y_i )
\end{equation}
when $i$ goes to infinity, that we denote by $l_+$ (this constant depends on the initial hyperplane that we choose, i.e.
$P_{C_+}$, $C_+ =(0,Y_+ )$), then we can conclude the proof using the same argument as in the proof of theorem \ref{main},
as long as $(x_- +l_+)_i \not= 0$ for $i=1,\dots ,n_-$.\\

Indeed, we have by taking the limit in (\ref{exten}):
\begin{equation}
\label {secondexten}
\left .
\begin{array}{l}
\di\sum_{w\in \N^{n_+}} a_{0,w} Y_+^w =e ,\\
\di\sum_{w\in \N^{n_+}} a_{\nu_j ,w} (x_- +C_+ )^{\nu_j} Y_+^w =0 .
\end{array}
\right .
\end{equation}
Choosing $Y_+$ of the form
\begin{equation}
Y_+ (n)=\delta^n C ,
\end{equation}
with $\delta$ satisfying a multiplicative non resonance condition and such that $(x_- +l_+)_i \not= 0$ for $i=1,\dots ,n_-$,
we conclude that $I$ is trivial on $U$.\\

We note that the condition for $x_- +l_+$ to be admissible can be easily satisfied, indeed $l_+$ varies continuously with
$Y_+$. This follows from the fact that $W_+ (p)$ contains a graph portion $\Lambda$ near the homoclinic point.
As a consequence, we can choose $Y_+$ such that $(x_- +l_+ )_i \not= 0$ for $i=1,\dots ,n_-$.\\

This is not so clear for the condition given by equation (\ref{condi}). We study this condition in the following section.

\subsection{Study of the remainder condition}

We have to find conditions under which condition (\ref{condi}) can be satisfied.\\

We denote by $f_{\rm lin} (x,y) =(\lambda_- x ,\lambda_+ y)$. Then, we have
\begin{equation}
r_i (x,y) =\di\sum_{j=1}^i f_{\rm lin}^{j-1} \circ r \circ f^{i-j} (x,y) .
\end{equation}
we deduce, denoting $r_i (x,y)=(r_{i,-} (x,y) ,r_{i,+} (x,y))$, that
\begin{equation}
r_{i,-} (x,y)=\di\sum_{j=1}^i \lambda_-^{j-1} r_- (f^{i-j} (x,y)) .
\end{equation}
The quantity (\ref{condi}) can be written
\begin{equation}
\lambda_-^{-i} r_{-,i} (x_i ,y_i ) =
\di\sum_{j=1}^i \lambda_-^{j-i-1} r_- (f^{i-j} (x_i , y_i )) .
\end{equation}

As $f^i (x_i ,y_i ) =(X_i ,Y_+ )$ and $(X_i ,Y_+ )\rightarrow (0,Y_+)$ when $i$ goes to infinity, the quantity
$\mid r_- (f^{i-j} (x_i ,y_i )) \mid$ is always bounded. As $W^+ (p)$ contains a graph portion $\Lambda$, we have
$x_i =\Lambda (y_i)$ and $y_i$ depends continuously on $Y_+$. Then, the quantity $\lambda_-^{-i} r_{-,i} (x_i ,y_i )$
is bounded for all $i$ and admits a limit $l_+$ when $i$ goes to infinity which depends continuously on $Y_+$.\\

We then are lead to introduce the following class of diffeomorphisms:

\begin{defi}[Normal form]
Let $f$ be an analytic diffeomorphism of $\R^n$ possessing an hyperbolic point $p$ with stable and unstable manifolds of dimension
$n_-$ and $n_+$ respectively, such that $n_- +n_+ =n$. A normal form for $f$ is an analytic coordinates system $(x,y)$
defined on an open neighbourhood $U$ of $p$ such that $f$ takes the form
\begin{equation}
f(x,y)=(\lambda_- x ,\lambda_+ y) +(r_- (x,y), r_+ (x,y) ),
\end{equation}
with $\partial_y r_- =\partial_x r_+ =0$.
\end{defi}

\subsection{Main result}

The previous discussion leads to the following version of theorem \ref{main}:

\begin{thm}
\label{main2}
Let $f$ be an analytic diffeomorphism of $\R^n$, which possesses a hyperbolic point $p$, with a stable manifold denoted
by $W^- (p)$ and an unstable manifold $W^+ (p)$ of dimension $n_-$ and $n_+$ such that $n_- +n_+ =n$. We denote by
$\lambda^-_i$, $i=1,\dots ,n_-$ and $\lambda^+_i$, $i=1,\dots ,n_+$ the eigenvalues of $Df(p)$ associated to $W^- (p)$ and
$W^+ (p)$. We assume that:\\

i) $W^+ (p)$ and $W^- (p)$ do not coincide and intersect along (at least) one homoclinic orbit $\gamma$,

ii) $f$ is in normal form on an open neighbourhood $U$ of $p$.

iii) We denote by $h=(h^- ,h^+ )$ the first intersection of $\gamma$ with $U$.
The spectrum $(\lambda^{\sigma})$ satisfies a multiplicative non resonance condition and $h_i^{\sigma} \not= 0$ for
$i=1,\dots ,n_{\sigma}$, for $\sigma =+$ or $-$.

iv) $W^+ (p)$ contains a graph portion $\Lambda$ associated with $h$ in a small neighbourhood $V\subset U$.\\

Then, the dynamical system defined by $f$ does not possess a non trivial analytic first integral.
\end{thm}

The main advantage of this theorem is that it covers some symplectic cases. However, we can probably proved a more general
theorem by restricting our attention to symplectic diffeomorphisms. Indeed, in this case we have more information on the geometry
of the stable and unstable manifolds as well as on the remainder of the normal form. We hope that this method can be used to solve the
cases which are not considered in \cite{cdr} leading to a complete proof of the Birkhoff conjecture. This will be studied in a forthcoming paper.

\end{document}